\documentclass[11pt]{article}
\usepackage[a4paper]{anysize}\marginsize{2.cm}{2.0cm}{2.0cm}{2.5cm}
\pdfpagewidth=\paperwidth \pdfpageheight=\paperheight
\usepackage{amsfonts,amssymb,amsthm,amsmath,eucal,tabu,url}
\usepackage{pgf}
 \usepackage{array}
 \usepackage{pstricks}
 \usepackage{pstricks-add}
 \usepackage{pgf,tikz}
 \usetikzlibrary{automata}
 \usetikzlibrary{arrows}
 \usepackage{indentfirst}
\theoremstyle{plain}
\newtheorem{thm}{Theorem}[section]
\newtheorem{theorem}[thm]{Theorem}

\newtheorem{proposition}[thm]{Proposition}
\newtheorem{corollary}[thm]{Corollary}

\theoremstyle{definition}

\newtheorem{remark}[thm]{Remark}

\newtheorem{question}[thm]{Question}
\newtheorem{problem}[thm]{Problem}

\newtheorem{thevarthm}[thm]{\varthmname}

\newenvironment{varthm*}[1]{\trivlist\item[]{\bf #1.}\it}{\endtrivlist}

\renewcommand\geq{\geqslant}

\renewcommand\leq{\leqslant}

\newcommand\be{\begin{eqnarray*}}
\newcommand\ee{\end{eqnarray*}}

\newcommand\newop[2]{\def#1{\mathop{\rm #2}\nolimits}}
\newop\log{log}
\newop\ord{ord}
\newop\Gal{Gal}
\newop\SL{SL}
\newop\Bl{Bl}
\newop\mult{mult}
\newop\mass{mass}
\newop\div{div}
\newop\codim{codim}
\newop\sing{sing}
\newop\vdim{vdim}
\newop\edim{edim}
\newop\Ass{Ass}
\newop\size{size}
\newop\reg{reg}
\newop\satdeg{satdeg}
\newop\supp{supp}
\newop\Neg{Neg}
\newop\Nef{Nef}
\newop\Nefh{Nef_H}
\newop\Eff{Eff}
\newop\Zar{Zar}
\newop\MB{MB}
\newop\MBxC{MB\mathit{(x,C)}}
\newop\NnB{NnB}
\newop\Bigg{Big}
\newop\Effbar{\overline{\Eff}}

\def\keywordname{{\bfseries Keywords}}%
\def\keywords#1{\par\addvspace\medskipamount{\rightskip=0pt plus1cm
\def\and{\ifhmode\unskip\nobreak\fi\ $\cdot$
}\noindent\keywordname\enspace\ignorespaces#1\par}}
\def\subclassname{{\bfseries Mathematics Subject Classification
(2000)}\enspace}
\def\subclass#1{\par\addvspace\medskipamount{\rightskip=0pt plus1cm
\def\and{\ifhmode\unskip\nobreak\fi\ $\cdot$
}\noindent\subclassname\ignorespaces#1\par}}

\begin{document}
\title{ON THE SYLVESTER-GALLAI AND THE ORCHARD PROBLEM FOR PSEUDOLINE ARRANGEMENTS}
\author{J\"urgen Bokowski and Piotr Pokora}
%
\date{\today}
\maketitle
\thispagestyle{empty}
\begin{abstract}                                            
We study a non-trivial extreme case of the orchard problem for $12$ pseudolines and we provide a complete classification of pseudoline arrangements having $19$ triple points and $9$ double points. We have also classified those that can be realized with straight lines. They include new examples different from the known example of B\"or\"oczky. 
Since Melchior's inequality also holds for arrangements of pseudolines, we are able to deduce that some combinatorial point-line configurations cannot be realized using pseudolines. In particular, this gives a negative answer to one of Gr\"unbaum's problems. We formulate some open problems which involve our new examples of line arrangements.
\subclass{52C30, 32S22}
\end{abstract}
%
\section{The Sylvester-Gallai problem}
We begin with a few definitions.
A {\bf line arrangement} in the real projective plane $\mathbb{P}^{2}_{\mathbb{R}}$ is a finite set of lines in $\mathbb{P}^{2}_{\mathbb{R}}$.
A {\bf pseudoline} in $\mathbb{P}^{2}_{\mathbb{R}}$ is a simple closed curve such that its removal does not cut $\mathbb{P}^{2}_{\mathbb{R}}$ in two 
connected components. A {\bf pseudoline arrangement} is a set of pseudolines in $\mathbb{P}^{2}_{\mathbb{R}}$ such that every pair
of pseudolines has precisely one point in common where the two curves intersect each other. A first book about
pseudoline arrangements was written by Gr\"unbaum \cite{Gr pseudolines}.
It has turned out later that pseudoline arrangements are isomorphic to reorientation classes of oriented matroids 
in rank 3. This implies a close connection of our investigation to the theory of oriented matroids. The interested
reader can find more about this relation in \cite{Bj} and \cite{Bo}.
For a given line arrangement, or for a given pseudoline arrangement, we count the number of points that are incident
with precisely $k$ lines or $k$ pseudolines, respectively, with $k \geq 2$, and we denote this number of the arrangement by $t_k$.
We call a point incident with precisely $r$ lines or with precisely $r$ pseudolines an {\bf $r$-point}. We use also
the notion {\bf double point} for $r=2$, {\bf triple point} for $r=3$, and {\bf quadruple point} for $r=4$.
We speak of an {\bf essential line (pseudoline) arrangement} when all lines (pseudolines) do not intersect at one point.
Our article can also be seen in the spirit of Gr\"unbaum's book about point-line configurations, see \cite{Gruenbaum}.
The reader will get some benefit for understanding our paper when she/he has a look at this book. 
For a {\bf point-line configuration} we have not only a set of lines but also a set of points 
together with an incidence relation between the set of points and the set of lines. 
It is clear that the lines can be replaced with pseudolines and we arrive
at a {\bf point-pseudoline arrangement}. Moreover, when we forget about any underlying geometric set of points, lines, 
or pseudolines, we arrive at an {\bf abstract point-line configuration}.  
Point-line configurations with $n$ lines and $n$ points in which 
$(\bullet)$ the lines are incident with precisely $k$ points and
$(\bullet \bullet)$ the points are incident with precisely $k$ lines 
have been called $(n_k)${\bf -configurations}. We refer to them later on.  
We use the point-line duality of the projective plane, i.e., a line arrangement defines via duality a point configuration, and vice versa. 

The solution of the famous problem due to Sylvester \cite{Sylvester} says that for every finite configuration of points in
the real projective plane there exists at least one \textbf{ordinary line}, i.e., a line passing through exactly two points from the configuration, provided that not all points lie on a line. The dual version of this problem 
tells us that every arrangement of lines in the real projective plane, not all intersecting at one point, contains at least one double intersection point, 
           i.e., $t_{2}\geq 1$. 
Sylvester's problem was solved by Gallai \cite{Gallai}. It is worth pointing out that Melchior \cite{M41} has shown in the dual situation that for essential line arrangements of at least $3$ lines one has $t_{2} \geq 3$. It was natural to ask what is the maximal possible number of ordinary lines for configurations of points or, via the famous orchard problem \cite{Jackson}, what is the maximal possible number of triple intersection points for line arrangements. Quite recently, Green and Tao \cite{GreTao13} have shown the so-called \emph{Dirac-Motzkin} conjecture which provides also the upper bound for the number of triple intersection points for arrangements of $n \gg 0$ lines defined over the reals, namely
%
%
$$t_{3} \leq 1+ \left\lfloor \frac{n(n-3)}{6} \right\rfloor.$$
It was well-known several years before the proof of Green and Tao that their upper-bound can be obtained using the so-called \emph{B\"or\"oczky family of line arrangements} \cite{FP84}. 
Sylvester's problem provides also a lot of geometrical constraints, for instance it implies that the famous
\emph{dual-Hesse arrangement} of $9$ lines and $12$ triple points (see for instance \cite{AD}) cannot be realized as a straight-line arrangement in the real projective plane since it does not contain any double point. It means also, via duality, that the \emph{Hesse arrangement} of $12$ lines, $12$ double points, and $9$ quadruple points cannot be realized as a straight-line arrangement in the real projective plane. 
%

An analogon to Sylvester's problem can be formulated for pseudoline arrangements in which \emph{not all pseudolines intersect in a common point}.
\begin{problem} 
\label{pseudolinesorchard}
Let $\mathcal{L}$ be a {\it pseudoline arrangement} with $n\geq 3$ pseudolines not intersecting in a common point. Is it true that $t_{2} \geq 1 $?
\end{problem}
Probably questions around Sylvester's problem and the orchard problem for pseudolines were considered for the first time in the paper due to Burr, Gr\"unbaum, and Sloane \cite{BGS74}. However, the authors did not know whether there exists a duality result for configurations of points and arrangements of pseudolines. In 1980 Goodman \cite[Theorem~2]{Goodman} has shown the following result.
\begin{theorem}(Duality \cite{Goodman})
\label{Goodman}
If, in $\mathbb{P}^{2}_{\mathbb{R}}$, $\mathcal{L}$ is an arrangement of pseudolines and $\mathcal{P}$ a configuration of points, and if $I$ is the set of all true statements of the form ``$P\in \mathcal{P}$ is incident to $L\in \mathcal{L}$", then there is a configuration $\mathcal{L}'$ of points and an arrangement $\mathcal{P}'$ of pseudolines, such that the set of all incidences holding between members of $\mathcal{L}'$ and members of $\mathcal{P}'$ is the dual $I'$ of $I$.
\end{theorem}
This result allows to establish some bounds in the context of the orchard problem \cite{BGS74}.

On the other hand, let us point out here explicitly that Problem \ref{pseudolinesorchard} has a positive answer. Indeed, using the same proof as in the case of Melchior's inequality, one can show the following result. In our formulation by $p_{j}$ we mean the number of regions bounded by precisely $j$ sections of pseudolines.
\begin{theorem}(Melchior \cite{M41}) Let $\mathcal{L}$ be an essential arrangement of $n \geq 3$ pseudolines. Then
$$\sum_{r \geq 2}(3-r)t_{r} = 3 + \sum_{j\geq 3}(j-3)p_{j}.$$
\end{theorem}
\begin{corollary}
For an essential arrangement $\mathcal{L}$ of $n \geq 3$ pseudolines one has
$$t_{2} \geq 3.$$
\end{corollary}
\begin{remark}
The fact that every configuration of pseudolines has at least one double point follows from a result by Kelly and Rottenberg \cite{KR}. However, our aim was to provide an explicit inequality which is more adequate for our purposes.
\end{remark}
This quite natural result has some significant geometrical consequences. Before we formulate some corollaries, it is worth pointing out that for arrangements of $n$ pseudolines we have the same combinatorial equality as in the case of straight lines, namely
$${n \choose 2}= \sum_{r \geq 2} {r \choose 2}t_{r}.$$
On the left hand side of the equation, we have the number of pairwise intersections, and on the right hand side, we have the sum over all $r$-points in the arrangement.
Using this fact we can derive the following result.
%
%
\begin{corollary}
There does not exist an arrangement of $n=9$ pseudolines with $12$ triple points.
\end{corollary}
\begin{proof}
Using the above combinatorial equality, observe that if $n=9$ and $t_{3}=12$, then $t_{r} = 0$ for $r\neq 3$. In particular, $t_{2} = 0$.
\end{proof}
Now by Theorem \ref{Goodman}, one has the following.
\begin{corollary}
There does not exist an arrangement of $12$ pseudolines intersecting at $9$ quadruple points and $12$ double points.
\end{corollary}
At last, let us recall the following question which was formulated by Gr\"unbaum in his book \cite{Gruenbaum}. Before we proceed further, let us also recall that by a \textbf{geometric realization} of a given abstract point-line configuration we mean a realization with straight lines, and by a \textbf{topological realization} we mean a realization with pseudolines. 
\begin{problem}(\cite[Page~254, Problem 1]{Gruenbaum})
Decide whether any abstract point-line configuration of $26$ triple points and $13$ lines can be realized geometrically or
topologically.
\end{problem}
Now we disprove partially the above problem of Gr\"unbaum.  
\begin{proposition}
Abstract point-line configurations of $26$ triple points and $13$ lines are not realizable with pseudolines.
\end{proposition}
\begin{proof}
Using the combinatorial equality for pseudoline arrangements, we see that for $13$ pseudolines $26$ triple points is the maximal possible number of intersection points, which means that $t_{r}=0$ for $r\neq 3$. This fact, combined with
Melchior's inequality, completes the proof.
\end{proof}
Using exactly the same argument, one can show that:
\begin{enumerate}
\item \emph{Klein's arrangement} \cite{Klein} consisting of $21$ lines and $t_{3} = 28$, $t_{4}=21$;
\item \emph{Wiman's arrangement} \cite{Wiman96} consisting of $45$ lines and $t_{3}=120, t_{4}=45, t_{5}=36$;
\item \emph{Fermat's family of line arrangements} \cite[Example II.6]{GU} for $n\geq 4$ consisting of $3n$ lines and $t_{n}=3$, $t_{3}=n^{2}$;
\end{enumerate}
cannot be constructed as pseudoline arrangements. Our considerations here lead to the following very interesting question which seems to be not formulated explicitly in the literature. 
\begin{problem}
Let $\mathcal{L}$ be an arrangement of $n\geq 4$ lines in the complex projective plane
$\mathbb{P}^{2}_{\mathbb{C}}$ 
such that $t_{2} = 0$. Then $\mathcal{L}$ is isomorphic to either the dual Hesse, Klein's, Wiman's arrangement, or to one of Fermat's arrangements of lines.  
\end{problem}
On the other hand, as the referee kindly pointed out to us, it might be interesting to ask whether the dual Hesse, Klein's, Wiman's arrangement, or to one of Fermat's arrangements of lines, can be realized as line arrangements in the three dimensional real projective space -- this is a natural intermediate step between the real projective plane and the complex projective plane.
%
%
\section{The orchard problem for pseudoline arrangements}
%
\subsection{B\"or\"oczky's arrangement of 12 lines}
%
Let us here recall the main construction due to B\"or\"oczky which motivated our research -- for more details please consult \cite{FP84}. 

We start with a regular $n$-gon inscribed in a circle $O$. We denote vertices of this $n$-gon traced in clockwise order by $P_{0},\ldots, P_{n-1}$. For simplicity, we assume here that $n$ is even --
we are only interested in the case where $n$ is even and thus we omit the odd case.
Then we construct the first line by joining $P_{0}$ with $P_{n/2}$. In the next step we join $P_{n/2 - 2}$ with $P_{1}$, and we continue this procedure until we obtain exactly $n$ lines -- after this moment our construction repeats and it does not provide new distinct lines (since we consider indices modulo $n$). Of course, it may happen that $P_{n/2 - 2i}$ and $P_{i}$ coincide 
-- then we draw the tangent line at $P_{i}$ to $O$. We obtain the arrangement $\mathcal{B}_{n}$ which consists of $n$ lines, $n-3 +\varepsilon$ double points, and $1 + \left\lfloor \frac{n(n-3)}{6} \right\rfloor$ triple points, where
$\varepsilon$ is equal to $0$ if $n = 0 \, {\rm mod}(3)$, or $2$, otherwise.
%
%
\begin{figure}[h]
\centering
\definecolor{uuuuuu}{rgb}{0.27,0.27,0.27}
\begin{tikzpicture}[line cap=round,line join=round,x=1.0cm,y=1.0cm,scale=0.7]
\clip(-0.82,-3.96) rectangle (8.3,5.22);
\draw [domain=-4.36:18.32] plot(\x,{(-2.33-0.03*\x)/-2.89});
\draw [domain=-4.36:18.32] plot(\x,{(-5.21--0.51*\x)/-1.98});
\draw [domain=-4.36:18.32] plot(\x,{(--7.45-1.43*\x)/1.46});
\draw [domain=-4.36:18.32] plot(\x,{(--8.44-2.49*\x)/1.47});
\draw [domain=-4.36:18.32] plot(\x,{(--4.08-1.97*\x)/0.55});
\draw [domain=-4.36:18.32] plot(\x,{(-3.23--1.98*\x)/0.51});
\draw [domain=-4.36:18.32] plot(\x,{(-6.11--2.51*\x)/1.42});
\draw [domain=-4.36:18.32] plot(\x,{(-5.12--1.46*\x)/1.43});
\draw [domain=-4.36:18.32] plot(\x,{(--2.03-0.55*\x)/-1.97});
\draw [domain=-4.36:18.32] plot(\x,{(--5.19-0.71*\x)/1.26});
\draw [domain=-4.36:18.32] plot(\x,{(--3.18-0.73*\x)/-1.24});
\draw [domain=1:2] plot(\x,{(-2.12--1.44*\x)/-0.01});
\draw(2.9,0.83) circle (1.44cm);
\begin{scriptsize}
\fill [color=uuuuuu] (2.19,-0.42) circle (2.0pt);
\fill [color=uuuuuu] (4.35,0.85) circle (2.0pt);
\fill [color=uuuuuu] (2.17,2.08) circle (2.0pt);
\fill [color=uuuuuu] (0.9,4.23) circle (2.0pt);
\fill [color=uuuuuu] (1.44,3.32) circle (2.0pt);
\fill [color=uuuuuu] (1.45,2.26) circle (2.0pt);
\fill [color=uuuuuu] (2.36,2.8) circle (2.0pt);
\fill [color=uuuuuu] (1.85,0.82) circle (2.0pt);
\fill [color=uuuuuu] (2.9,0.83) circle (2.0pt);
\fill [color=uuuuuu] (1.48,-1.68) circle (2.0pt);
\fill [color=uuuuuu] (0.96,-2.6) circle (2.0pt);
\fill [color=uuuuuu] (2.39,-1.14) circle (2.0pt);
\fill [color=uuuuuu] (1.47,-0.62) circle (2.0pt);
\fill [color=uuuuuu] (3.44,-0.08) circle (2.0pt);
\fill [color=uuuuuu] (4.88,0.32) circle (2.0pt);
\fill [color=uuuuuu] (4.87,1.38) circle (2.0pt);
\fill [color=uuuuuu] (5.79,0.86) circle (2.0pt);
\fill [color=uuuuuu] (6.85,0.87) circle (2.0pt);
\fill [color=uuuuuu] (3.42,1.75) circle (2.0pt);
\fill [color=uuuuuu] (3.61,2.09) circle (2.0pt);
\fill [color=uuuuuu] (3.64,-0.41) circle (2.0pt);
\fill [color=uuuuuu] (1.46,0.82) circle (2.0pt);
\end{scriptsize}
\end{tikzpicture}
\caption{B\"or\"oczky's arrangement of $12$ lines.}
\end{figure}
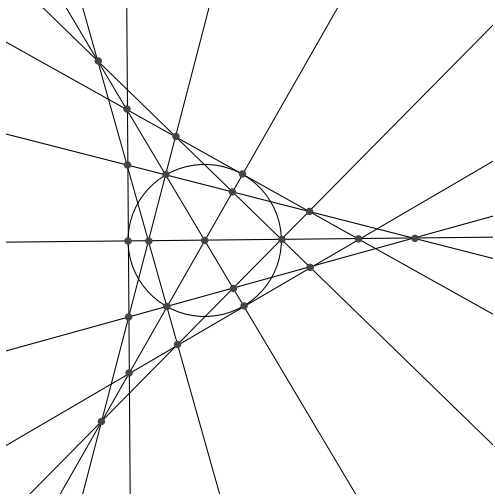
\subsection{Sweeping argument}
Now we explain how we determine pseudoline arrangements with 12 pseudolines in which $19$ triples of them intersect in a point and
in which the remaining pairs intersect pairwise.
We can assume that one pseudoline $P_{1}$ is incident with $5$ triple
points since $3\cdot t_{3}/n = 3\cdot 19/12 > 4$. We use an additional {\it sweeping} pseudoline $P$ that is at the beginning equal to pseudoline $P_1$
and that has otherwise all the time precisely one additional point on $P_{1}$ between the intersection of $P_{12}$ with $P_{1}$ and the intersection of $P_{2}$ with 
$P_{1}$. 
We sweep $P$ through the projective plane until it finally coincides with $P_1$ again. We move $P$ across all
intersection points of two or three pseudolines. In other words, we imagine all possible crossings that can occur
when we have an arrangement with 12 pseudolines in which $19$ triples occur and the remaining crossings occur
pairwise.  
This  sweeping process has been used heavily in \cite{BP2} and we refer the reader for details
of this sweeping process to this article.
On this sweeping pseudoline $P$, we consider the sequence of intersections along $P$ with the remaining pseudolines
$P_{i}$ for $2 \leq i \leq 12$. An intersection of just two pseudolines in the arrangement leads to an interchange of just two elements in the sequence of intersections on $P$. An intersection of three pseudolines in the arrangement leads to reversing the sequence of a triple of adjacent indices. 
We have at the beginning on $P$ the sequence of pseudoline indices $(2,3,4,5,6,7,8,9,10,11,12)$ and immediately
afterwards when $P$ is incident with line $P_1$, we have triples and a single crossing point
$((1,2,3),(1,4,5),(1,6,7),(1,8,9),(1,10,11),12)$ or the new sequence of pseudoline
indices $(3,2,5,4,7,6,9,8,11,10,12)$. Because of the symmetry, we
consider only one position of $P_{12}$. Now we change either a pair of
adjacent indices in ascending order or we look for a triple of adjacent indices
again in ascending order for which we reverse their order. We do such steps
for changing the sequence of indices over and over again in all possible ways
until we reach the sequence of indices $(12,11,10,9,8,7,6,5,4,3,2)$. 
During the process we collect the triples and we keep those examples that have
precisely $19$ triples.

We delete those elements for which we have already an isomorphic example,
and we decide which of them are realizable with straight lines. The functional
programming language Haskell \cite{Haskell} (our code is available on demand) was used for this algorithm,
and Cinderella \cite{Cinderella} software  was helpful to investigate the outcome.
The sweeping process that was used in \cite{BP} was more involved because of the problem size.
However, the motivation for both investigations can be seen as an equal
attempt to work towards a solution of the long standing open problem to
get a complete classification of all geometric $(n_4)$-configurations. 

In Section 3, we provide a complete list of arrangements with $12$ pseudolines, $19$ triple points, and $9$
double points -- this list was obtained using machine-computations based on Haskell code. For all those that are realizable with straight lines, we provide in Section 4 corresponding pictures. We also explain why some of them
cannot be represented by lines.
%
\subsection{Problems and Applications}
The methods we have used to classify all arrangements with $12$ pseudolines and $19$ triple points are
useful in the context of classifying $(n_4)$-configurations in which long standing problems are still open. 
One variant of the sweeping method was decisive for investigations of the case $n=19$, see \cite{BP}.  
Apart from providing building blocks for this area, our investigated case is as well strongly related to the so-called {\it containment problems} in algebraic geometry that we are going to describe next. Here we focus on a special case of this problem which fits to our setting -- for a general introduction please consult a very recent 
%
%
survey~\cite{SS}.

Let $\mathcal{P} = \{P_{1}, \ldots, P_{s}\}$ be a finite set of mutually distinct points in the projective plane and denote by
%
$I(P_{i})$ the \textbf{associated radical ideal} of $P_{i}$, that is the ideal of rational functions vanishing on $P_{i}$. We define the radical ideal $I = I(\mathcal{P})$ of $\mathcal{P}$ by

$$I = I(P_{1}) \cap \ldots \cap I(P_{s}).$$
The \textbf{$m$-th symbolic power of $I$} is defined as
$$I^{(m)} = I^{m}(P_{1}) \cap \ldots \cap I^{m}(P_{s}),$$
which means that the $m$-th symbolic power of $I$ can be viewed as the set of forms vanishing along points $P_{i}$ with multiplicity $\geq m$. Another famous result due to Ein, Lazarsfeld, and Smith \cite{ELS} in characteristic $0$, and Hochster and Huneke \cite{HH} in positive characteristic, tells us that one has the the following containment:
$$I^{(2k)} \subset I^{k}.$$ 
Few years later, Huneke asked whether one can improve the above containment.
\begin{problem}(\cite[Problem 0.4]{Huneke}) Does the containment
\begin{equation}
\label{cont}
I^{(3)} \subset I^{2}
\end{equation}
hold?
\end{problem}
It turned out that in general (\ref{cont}) does not hold and most counter-examples are based on radical ideals of points which are given by intersection points of line arrangements. In particular, in \cite{CG13} the very first 
counter-example to (\ref{cont}) in the real projective plane is provided and it is given by the radical ideal of triple points of B\"or\"oczky's arrangement of $12$ lines. In order to understand better the containment problem, one can formulate the following question.
\begin{question} Using singular intersection points of line arrangements from Section 4, decide whether taking radical ideals of triple points one always has the containment
$$I^{(3)} \subset I^{2}.$$
\end{question}
This question is of interest because the mentioned configurations have the same combinatorics, 
i.e., the number of lines and types of singular points, 
thus the main problem is whether the containment problem is combinatorial in nature.
%
We want to verify whether the following is true:
If $\mathcal{L}$ and $\mathcal{L}'$ are two line configurations with the same combinatorics (with the same number of lines $n$ and with the same combinatorial structure $(t_{2}, \ldots , t_{n})$) 
and if the containment (\ref{cont}) does not hold for the radical ideal of a certain subset of singular points of $\mathcal{L}$, 
then (\ref{cont}) does not hold 
for the radical ideal of corresponding singular points of $\mathcal{L}'$.
%
%
\section{Pseudoline arrangements} 
\begin{figure}[ht]
\centering 
\definecolor{sqsqsq}{rgb}{0.12549019607843137,0.12549019607843137,0.12549019607843137}
\definecolor{qqffqq}{rgb}{0.,1.,0.}
\definecolor{yqqqqq}{rgb}{0.5019607843137255,0.,0.}
\definecolor{ffqqqq}{rgb}{1.,0.,0.}
\definecolor{qqffff}{rgb}{0.,1.,1.}
\definecolor{xfqqff}{rgb}{0.4980392156862745,0.,1.}
\definecolor{qqwuqq}{rgb}{0.,0.39215686274509803,0.}
\definecolor{cqcqcq}{rgb}{0.7529411764705882,0.7529411764705882,0.7529411764705882}
\definecolor{ffxfqq}{rgb}{1.,0.4980392156862745,0.}
\definecolor{qqqqff}{rgb}{0.,0.,1.}
\definecolor{ffffqq}{rgb}{1.,1.,0.}
\begin{tikzpicture}[line cap=round,line join=round,>=triangle 45,x=1.0cm,y=1.0cm,scale=0.9]
\clip(-8.71178591366661,-5.522103959214514) rectangle (16.45454955189041,6.87212732111986);
\draw [line width=2.pt] (1.3904744381092704,0.44141752156722525) circle (5.679341127680998cm);
\draw [line width=2.pt,color=ffffqq] (-3.796387676265502,-1.871886836063335)-- (-2.544900105890475,-1.8340644075052521);
\draw [line width=2.pt,color=ffffqq] (-2.544900105890475,-1.8340644075052521)-- (1.3904744381092704,0.44141752156722525);
\draw [line width=2.pt,color=ffffqq] (1.3904744381092704,0.44141752156722525)-- (2.7931386610544315,1.1336806833224882);
\draw [line width=2.pt,color=ffffqq] (2.7931386610544315,1.1336806833224882)-- (4.838186822309607,2.7630686492005925);
\draw [line width=2.pt,color=ffffqq] (4.838186822309607,2.7630686492005925)-- (6.173334265099898,3.5039601775849425);
\draw [line width=2.pt] (7.049879223445249,0.9168679609245864)-- (4.838186822309607,2.7630686492005925);
\draw [line width=2.pt,color=qqqqff] (4.309837061621481,5.312991960802325)-- (4.838186822309607,2.7630686492005925);
\draw [line width=2.pt,color=ffxfqq] (-4.282010799452786,0.1624428377526984)-- (-2.893092811703862,0.2857338847532707);
\draw [line width=2.pt,color=cqcqcq] (-2.893092811703862,0.2857338847532707)-- (-2.544900105890475,-1.8340644075052521);
\draw [line width=2.pt,color=ffxfqq] (-2.893092811703862,0.2857338847532707)-- (1.8081957584220902,2.685295534958135);
\draw [line width=2.pt] (-4.282010799452786,0.1624428377526984)-- (-1.9952667896893945,2.99583835861175);
\draw [line width=2.pt] (-1.9952667896893945,2.99583835861175)-- (0.9974866170254968,4.143060497852456);
\draw [line width=2.pt,color=qqwuqq] (1.8081957584220902,2.685295534958135)-- (4.309837061621481,5.312991960802325);
\draw [line width=2.pt,color=cqcqcq] (-1.9952667896893945,2.99583835861175)-- (-2.893092811703862,0.2857338847532707);
\draw [line width=2.pt,color=cqcqcq] (-1.6268540015754405,5.252929210385896)-- (-1.9952667896893945,2.99583835861175);
\draw [line width=2.pt,color=cqcqcq] (-2.544900105890475,-1.8340644075052521)-- (-0.4656372707017834,-3.2556909798593434);
\draw [line width=2.pt,color=cqcqcq] (-0.4656372707017834,-3.2556909798593434)-- (4.2814835235882915,-4.447036437737969);
\draw [line width=2.pt,color=xfqqff] (0.9974866170254968,4.143060497852456)-- (0.7827881349174717,2.1358256589867013);
\draw [line width=2.pt,color=xfqqff] (1.8352000799560861,6.103319564878537)-- (0.9974866170254968,4.143060497852456);
\draw [line width=2.pt,color=xfqqff] (-0.4656372707017834,-3.2556909798593434)-- (1.3904744381092708,-5.237923606113773);
\draw [line width=2.pt,color=qqqqff] (-0.4656372707017834,-3.2556909798593434)-- (1.8288070077796332,-2.4077441812901257);
\draw [line width=2.pt,color=qqqqff] (1.8288070077796332,-2.4077441812901257)-- (4.904692453569939,-0.927993885747766);
\draw [line width=2.pt,color=qqqqff] (4.904692453569939,-0.927993885747766)-- (4.838186822309607,2.7630686492005925);
\draw [line width=2.pt,color=qqffff] (-3.2642904826045185,3.6953504000177366)-- (-2.893092811703862,0.2857338847532707);
\draw [line width=2.pt,color=qqffff] (1.8288070077796332,-2.4077441812901257)-- (6.2262140249722435,-2.5369273699144834);
\draw [line width=2.pt,color=ffqqqq] (-3.2642904826045185,3.6953504000177366)-- (-1.9952667896893945,2.99583835861175);
\draw [line width=2.pt,color=ffqqqq] (2.7931386610544315,1.1336806833224882)-- (4.904692453569939,-0.927993885747766);
\draw [line width=2.pt,color=ffqqqq] (4.904692453569939,-0.927993885747766)-- (6.2262140249722435,-2.5369273699144834);
\draw [line width=2.pt,color=yqqqqq] (-1.6268540015754405,5.252929210385896)-- (0.9974866170254968,4.143060497852456);
\draw [line width=2.pt,color=ffxfqq] (1.8081957584220902,2.685295534958135)-- (3.781474506869203,1.921135259662711);
\draw [line width=2.pt,color=ffxfqq] (3.781474506869203,1.921135259662711)-- (7.049879223445249,0.9168679609245864);
\draw [line width=2.pt,color=yqqqqq] (0.9974866170254968,4.143060497852456)-- (1.8081957584220902,2.685295534958135);
\draw [line width=2.pt,color=yqqqqq] (1.8081957584220902,2.685295534958135)-- (2.7931386610544315,1.1336806833224882);
\draw [line width=2.pt,color=yqqqqq] (2.7931386610544315,1.1336806833224882)-- (1.8288070077796332,-2.4077441812901257);
\draw [line width=2.pt,color=yqqqqq] (1.8288070077796332,-2.4077441812901257)-- (4.2814835235882915,-4.447036437737969);
\draw [line width=2.pt,color=qqwuqq] (-2.544900105890475,-1.8340644075052521)-- (-2.3976543663386622,-3.790012001121384);
\draw [line width=2.pt,color=qqqqff] (-2.3976543663386622,-3.790012001121384)-- (-0.4656372707017834,-3.2556909798593434);
\draw [line width=2.pt,color=qqffqq] (3.781474506869203,1.921135259662711)-- (4.904692453569939,-0.927993885747766);
\draw [line width=2.pt,color=qqffqq] (4.904692453569939,-0.927993885747766)-- (1.3904744381092708,-5.237923606113773);
\draw [line width=2.pt,color=qqffqq] (1.8352000799560861,6.103319564878537)-- (3.001247491289483,3.9384639129719976);
\draw [line width=2.pt,color=qqffqq] (3.001247491289483,3.9384639129719976)-- (3.781474506869203,1.921135259662711);
\draw [line width=2.pt] (0.9974866170254968,4.143060497852456)-- (3.001247491289483,3.9384639129719976);
\draw [line width=2.pt] (3.001247491289483,3.9384639129719976)-- (4.838186822309607,2.7630686492005925);
\draw [line width=1.2pt,color=ffqqqq] (-1.9952667896893945,2.99583835861175)-- (0.7857406789373714,2.1634292505996386);
\draw [line width=1.2pt,color=ffqqqq] (0.7857406789373714,2.1634292505996386)-- (2.7931386610544315,1.1336806833224882);
\draw [line width=2.pt,color=xfqqff] (0.7857406789373714,2.1634292505996386)-- (0.34204888745570544,0.5692072675839739);
\draw [line width=2.pt,color=xfqqff] (0.34204888745570544,0.5692072675839739)-- (-0.4656372707017834,-3.2556909798593434);
\draw [line width=2.pt,color=qqwuqq] (1.8081957584220902,2.685295534958135)-- (0.34204888745570544,0.5692072675839739);
\draw [line width=2.pt,color=qqwuqq] (-2.544900105890475,-1.8340644075052521)-- (0.34204888745570544,0.5692072675839739);
\draw [line width=2.pt,color=qqffff] (-2.893092811703862,0.2857338847532707)-- (0.34204888745570544,0.5692072675839739);
\draw [line width=2.pt,color=qqffff] (1.8288070077796332,-2.4077441812901257)-- (0.34204888745570544,0.5692072675839739);
\begin{scriptsize}
\draw [fill=black] (7.049879223445249,0.9168679609245864) circle (2.5pt);
\draw[color=black] (7.551863912723392,1.13601662491633) node {(1,2,3)};
\draw [fill=black] (-2.544900105890475,-1.8340644075052521) circle (2.5pt);
\draw[color=black] (-1.547315636985512,-1.7811622008209067) node {(7,11,12)};
\draw [fill=black] (-3.796387676265502,-1.871886836063335) circle (2.5pt);
\draw[color=black] (-3.5651569490887325,-1.2596975920751314) node {$P_{12}$};
\draw [fill=black] (2.7931386610544315,1.1336806833224882) circle (2.5pt);
\draw[color=black] (3.5161812885169508,1.2115912058939786) node {(4,6,12)};
\draw [fill=black] (4.838186822309607,2.7630686492005925) circle (2.5pt);
\draw[color=black] (5.632269555891114,2.828887238815659) node {(2,10,12)};
\draw [fill=black] (6.173334265099898,3.5039601775849425) circle (2.5pt);
\draw[color=black] (6.153734164636891,4.181672238315569) node {$P_{12}$};
\draw [fill=black] (4.309837061621481,5.312991960802325) circle (2.5pt);
\draw[color=black] (4.7102596679638005,5.8518704779216035) node {(1,10,11)};
\draw [fill=black] (-4.282010799452786,0.1624428377526984) circle (2.5pt);
\draw[color=black] (-3.648288988164146,-0.10340650311710745) node {(1,2,3)};
\draw [fill=black] (-2.893092811703862,0.2857338847532707) circle (2.5pt);
\draw[color=black] (-2.0007631228514047,1.1511315411118597) node {(3,5,7)};
\draw [fill=black] (1.8081957584220902,2.685295534958135) circle (2.5pt);
\draw[color=black] (2.6697459815672855,2.7986574064245997) node {(3,6,11)};
\draw [fill=black] (-1.9952667896893945,2.99583835861175) circle (2.5pt);
\draw[color=black] (-1.260132229270447,3.645092713374264) node {(2,4,7)};
\draw [fill=black] (0.9974866170254968,4.143060497852456) circle (2.5pt);
\draw[color=black] (0.5687726303886516,4.748481595647934) node {(2,6,9)};
\draw [fill=black] (-1.6268540015754405,5.252929210385896) circle (2.5pt);
\draw[color=black] (-1.6833498827452797,5.776295896943955) node {(1,6,7)};
\draw [fill=black] (-0.4656372707017834,-3.2556909798593434) circle (2.5pt);
\draw[color=black] (-0.9276040729687928,-2.4159886810331552) node {(7,9,10)};
\draw [fill=black] (4.2814835235882915,-4.447036437737969) circle (2.5pt);
\draw[color=black] (4.816064081332509,-4.592536613189436) node {(1,6,7)};
\draw [fill=black] (1.8352000799560861,6.103319564878537) circle (2.5pt);
\draw[color=black] (2.186068663310334,6.456467125742792) node {(1,8,9)};
\draw [fill=black] (1.3904744381092708,-5.237923606113773) circle (2.5pt);
\draw[color=black] (2.125608998528215,-4.89483493710003) node {(1,8,9)};
\draw [fill=black] (1.8288070077796332,-2.4077441812901257) circle (2.5pt);
\draw[color=black] (1.6268167640757336,-2.9752405802677546) node {(5,6,10)};
\draw [fill=black] (4.904692453569939,-0.927993885747766) circle (2.5pt);
\draw[color=black] (5.586924807304525,-0.8591523128935936) node {(4,8,10)};
\draw [fill=black] (-3.2642904826045185,3.6953504000177366) circle (2.5pt);
\draw[color=black] (-3.6785188205552055,4.143884947826744) node {(1,4,5)};
\draw [fill=black] (6.2262140249722435,-2.5369273699144834) circle (2.5pt);
\draw[color=black] (6.796118102946905,-2.5520229267929224) node {(1,4,5)};
\draw [fill=black] (3.781474506869203,1.921135259662711) circle (2.5pt);
\draw[color=black] (3.8335945286230753,2.3149800881676486) node {(3,8,12)};
\draw [fill=black] (-2.3976543663386622,-3.790012001121384) circle (2.5pt);
\draw[color=black] (-1.381051558834685,-3.882135551999538) node {(1,10,11)};
\draw [fill=black] (3.001247491289483,3.9384639129719976) circle (2.5pt);
\draw[color=black] (3.878939277209665,4.068310366849096) node {(2,8,11)};
\draw [fill=black] (0.7857406789373714,2.1634292505996386) circle (2.5pt);
\draw[color=black] (0.3269339712601757,2.723082825446951) node {(3,4,9)};
\draw [fill=sqsqsq] (0.34204888745570544,0.5692072675839739) circle (2.5pt);
\draw[color=sqsqsq] (1.0222201162545437,0.5163050608996111) node {(5,9,11)};
\end{scriptsize}
\end{tikzpicture}

\caption{The pseudoline arrangement $\mathcal{C}_{1}$.}
\label{pseudo}
\end{figure}
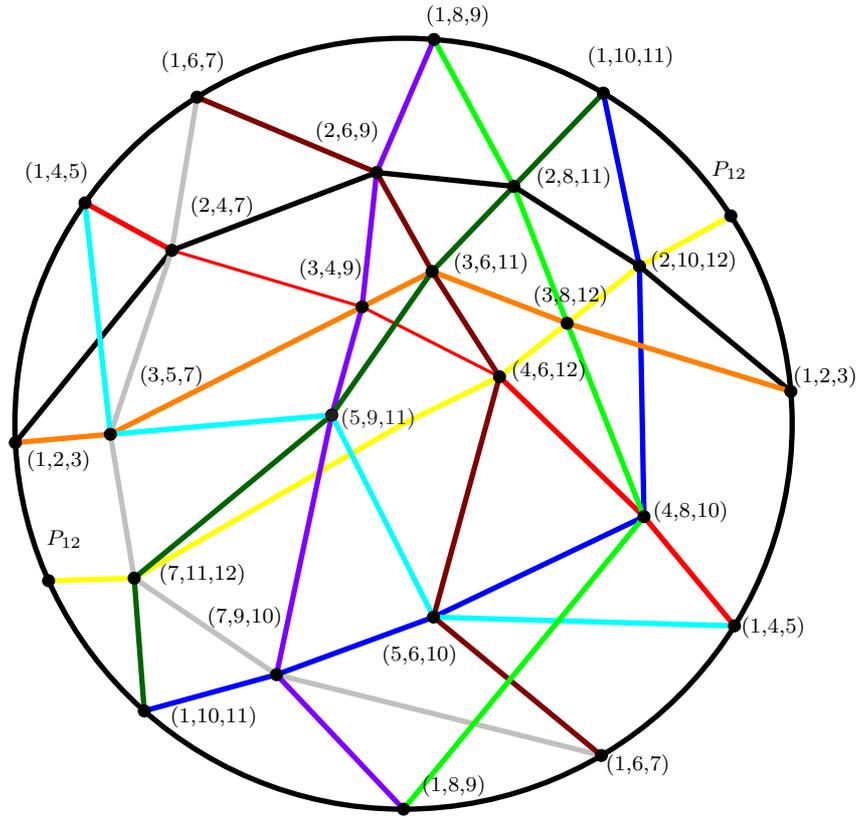

All arrangements with $12$ pseudolines, $19$ triple points, and $9$ double points can be reconstructed via
its triple points and this is one of the reasons why we provide only the triple lists. 
We can always start our sweeping algorithm with the line at infinity 
(we think of the circle model of the projective plane as in \cite{BP2})
with the triples $(1,2,3),(1,4,5),(1,6,7),(1,8,9),(1,10,11)$,
and we can insert successively all triples whereby one observes that no pair of pseudolines can meet twice. 
Finally, the list of $13$ non-isomorphic arrangements has appeared during our computer-based tests. 
For instance, in Figure \ref{pseudo} we have drawn the pseudoline arrangement $\mathcal{C}_{1}$ 
that cannot be realized with straight lines. Here you see the line at infinity as the boundary of the circle with identified
antipodal points.
We provide also a short explanation after the enumeration of all $13$ arrangements
why certain pseudoline arrangements are not realizable with straight lines. 
Proceeding along the same lines one can handle the remaining cases.
\small{
 \begin{multline*} \mathcal{C}_{1} = \{(1,2,3),(1,4,5),(1,6,7),(1,8,9),(1,10,11),(2,4,7),(2,6,9),(2,8,11),(2,10,12),(3,4,9), \\ (3,5,7),(3,6,11),(3,8,12),(4,6,12),(4,8,10),(5,6,10),(5,9,11),(7,9,10),(7,11,1)\},
\end{multline*}
\begin{multline*} \mathcal{C}_{2} = \{(1,2,3),(1,4,5),(1,6,7),(1,8,9),(1,10,11),(2,4,7),(2,6,9),(2,8,11),(2,10,12),(3,4,9), \\ (3,5,7),(3,6,11),(3,8,12),(4,6,12),(4,8,10),(5,9,11),(7,9,10),(7,11,12),(5,6,8)\},
\end{multline*}
 \begin{multline*} \mathcal{C}_{3} = \{(1,2,3),(1,4,5),(1,6,7),(1,8,9),(1,10,11),(2,4,7),(2,6,9),(2,8,11),(2,10,12),(3,5,7), \\ (3,6,11),(3,8,12),(4,6,12),(4,8,10),(4,9,11),(5,9,12),(7,9,10),(7,11,12),(5,6,8)\},
\end{multline*}
 \begin{multline*} \mathcal{C}_{4} = \{(1,2,3),(1,4,5),(1,6,7),(1,8,9),(1,10,11),(2,4,7),(2,6,9),(2,8,11),(2,10,12),(3,4,10), \\ (3,5,7),(3,6,12),(3,9,11),(4,6,11),(4,8,12),(5,6,8),(5,9,12),(7,11,12),(7,8,10)\},
\end{multline*}
 \begin{multline*} \mathcal{C}_{5} = \{(1,2,3),(1,4,5),(1,6,7),(1,8,9),(1,10,11),(2,4,7),(2,6,9),(2,8,11),(2,10,12),(3,4,10), \\ (3,5,7),(3,6,12),(3,9,11),(4,6,11),(4,8,12),(5,8,10),(5,9,12),(7,9,10),(7,11,12)\},
\end{multline*}
\begin{multline*} \mathcal{C}_{6} = \{(1,2,3),(1,4,5),(1,6,7),(1,8,9),(1,10,11),(2,4,7),(2,6,9),(2,8,11),(2,10,12),(3,4,9), \\ (3,5,7),(3,6,11),(3,8,12),(4,6,8),(4,11,12),(5,6,10),(5,9,12),(7,9,11),(7,8,10)\},
\end{multline*}
\begin{multline*} \mathcal{C}_{7} = \{(1,2,3),(1,4,5),(1,6,7),(1,8,9),(1,10,11),(2,4,7),(2,6,9),(2,10,12),(3,4,9),(3,5,7), \\ (3,6,8),(3,11,12),(4,6,11),(4,8,12),(5,6,10),(5,8,11),(5,9,12),(7,9,11),(7,8,10)\},
\end{multline*}
\begin{multline*} \mathcal{C}_{8} = \{(1,2,3),(1,4,5),(1,6,7),(1,8,9),(1,10,11),(2,4,7),(2,6,9),(2,8,11),(2,10,12),(3,4,10), \\ (3,5,7),(3,6,8),(3,9,11),(4,6,11),(4,8,12),(5,6,10),(5,9,12),(7,11,12),(7,8,10)\},
\end{multline*}
\begin{multline*} \mathcal{C}_{9} = \{(1,2,3),(1,4,5),(1,6,7),(1,8,9),(1,10,11),(2,4,7),(2,6,9),(2,8,11),(2,10,12),(3,4,10), \\ (3,6,8),(3,7,11),(3,9,12),(4,6,11),(4,8,12),(5,6,10),(5,7,9),(5,11,12),(7,8,10)\}.
\end{multline*}
\begin{multline*} \mathcal{C}_{10} = \{(1,2,3),(1,4,5),(1,6,7),(1,8,9),(1,10,11),(2,4,7),(2,6,9),(2,8,11),(3,4,9),(3,5,7), \\ (3,6,11),(3,8,10),(4,6,10),(4,8,12),(5,6,8),(5,9,11),(5,10,12),(7,11,12),(7,9,10)\}.
\end{multline*}
\begin{multline*} \mathcal{C}_{11} =  \{(1,2,3),(1,4,5),(1,6,7),(1,8,9),(1,10,11),(2,4,7),(2,6,9),(2,8,11),(2,10,12),(3,5,7), \\ (3,6,11),(3,8,12),(4,6,12),(4,8,10),(4,9,11),(5,6,10),(5,9,12),(7,9,10),(7,11,12)\},
\end{multline*}
\begin{multline*} \mathcal{C}_{12} = \{(1,2,3),(1,4,5),(1,6,7),(1,8,9),(1,10,11),(2,4,7),(2,6,9),(2,8,11),(2,10,12),(3,4,10), \\ (3,5,7),(3,6,12),(3,9,11),(4,6,11),(4,8,12),(5,6,10),(5,9,12),(7,11,12),(7,8,10)\},
\end{multline*}
\begin{multline*} \mathcal{C}_{13} = \{(1,2,3),(1,4,5),(1,6,7),(1,8,9),(1,10,11),(2,4,7),(2,6,9),(2,8,11),(2,10,12),(3,5,12), \\ (3,7,10),(3,9,11),(4,8,10),(4,11,12),(5,7,8),(5,9,10),(6,8,12),(7,9,12),(3,4,6)\}.
\end{multline*}
}
\noindent       
\begin{figure}[bh]
\begin{center}
\definecolor{qqwwtt}{rgb}{0.,0.4,0.2}
\definecolor{qqttzz}{rgb}{0.,0.2,0.6}
\definecolor{qqqqff}{rgb}{0.,0.,1.}
\definecolor{ffqqqq}{rgb}{1.,0.,0.}
\begin{tikzpicture}[line cap=round,line join=round,>=triangle 45,x=1.0cm,y=1.0cm,scale=0.8]
\clip(-10.015271276330946,-6.771725555060792) rectangle (9.43740812453375,2.9783949270351715);
\draw [domain=-10.015271276330946:9.43740812453375] plot(\x,{(-2.-2.*\x)/2.});
\draw [domain=-10.015271276330946:9.43740812453375] plot(\x,{(-6.--2.*\x)/2.});
\draw [domain=-10.015271276330946:9.43740812453375] plot(\x,{(-2.-0.*\x)/2.});
\draw (2.,-6.771725555060792) -- (2.,2.9783949270351715);
\draw (0.,-6.771725555060792) -- (0.,2.9783949270351715);
\draw [domain=-10.015271276330946:9.43740812453375] plot(\x,{(-6.-0.*\x)/2.});
\draw [color=qqttzz] (4.735482645648355,-6.771725555060792) -- (4.735482645648355,2.9783949270351715);
\draw [color=qqttzz,domain=-10.015271276330946:9.43740812453375] plot(\x,{(-6.470965291296711-1.*\x)/3.7354826456483554});
\draw [color=qqttzz,domain=-10.015271276330946:9.43740812453375] plot(\x,{(-0.011900013346615168-2.7354826456483554*\x)/-7.470965291296711});
\draw [dash pattern=on 1pt off 1pt,color=ffqqqq,domain=-10.015271276330946:9.43740812453375] plot(\x,{(--8.203262267387618-0.7322969760909079*\x)/-4.735482645648355});
\draw [->,color=ffqqqq] (4.735482645648355,-3.) -- (6.,-3.);
\draw [->,color=ffqqqq] (4.735482645648355,-1.) -- (6.,-1.);
\draw [->,color=ffqqqq] (4.735482645648355,1.7354826456483554) -- (5.59081743817813,2.5908174381781297);
\draw [->,color=ffqqqq] (-2.7354826456483554,-1.) -- (-3.77162705053996,-1.);
\draw [->,color=ffqqqq] (-8.197747451591562,-3.) -- (-9.490905040647489,-3.);
\draw [->,color=ffqqqq] (0.,-1.732296976090908) -- (0.,-1.302145772421593);
\begin{scriptsize}
\draw [fill=ffqqqq] (0.,-3.) circle (2.0pt);
\draw[color=ffqqqq] (0.4363822648426888,-3.2343342825930494) node {(1,2,3)};
\draw [fill=ffqqqq] (0.,-1.) circle (2.0pt);
\draw[color=ffqqqq] (0.44827265567451324,-0.6541194720871661) node {(1,4,5)};
\draw [fill=ffqqqq] (2.,-1.) circle (2.0pt);
\draw[color=ffqqqq] (1.5184078305387079,-0.6541194720871661) node {(3,5,7)};
\draw [fill=ffqqqq] (2.,-3.) circle (2.0pt);
\draw[color=ffqqqq] (1.5421886122023567,-3.2343342825930494) node {(2,4,7)};
\draw[color=black] (-3.850103630030002,2.5027792937621975) node {$P_{4}$};
\draw[color=black] (-2.827530018493105,-6.117754059310452) node {$P_{3}$};
\draw[color=black] (-8.951081296882665,-1.1832418641033495) node {$P_{5}$};
\draw[color=black] (1.7859416242547566,2.2649714771257106) node {$P_{7}$};
\draw[color=black] (-0.34243833464180845,2.253081086293886) node {$P_{1}$};
\draw[color=black] (-8.951081296882665,-2.669540718081393) node {$P_{2}$};
\draw [fill=qqqqff] (1.,-2.) circle (2.0pt);
\draw[color=qqqqff] (1.6492021296887762,-1.8312681644377766) node {(3,4,9)};
\draw [fill=qqqqff] (4.735482645648355,1.7354826456483554) circle (2.0pt);
\draw[color=qqqqff] (5.335223287554336,1.5456028318003379) node {(3,6,11)};
\draw[color=qqttzz] (4.282923698937878,2.6692447654077385) node {$P_{6}$};
\draw [fill=qqwwtt] (4.735482645648355,-3.) circle (2.0pt);
\draw[color=qqwwtt] (5.204428988404268,-2.6873763043291294) node {(2,6,9)};
\draw[color=qqttzz] (-9.272121849341923,1.099713175606925) node {$P_{9}$};
\draw [fill=qqqqff] (-2.7354826456483554,-1.) circle (2.0pt);
\draw[color=qqqqff] (-2.7740232597498955,-0.5471059546007471) node {(5,9,11)};
\draw[color=qqttzz] (-9.367244975996517,-3.6564431571228138) node {$P_{11}$};
\draw [fill=qqqqff] (0.,-1.732296976090908) circle (2.0pt);
\draw[color=qqqqff] (-0.5029586108714377,-2.045295199410615) node {(1,8,9)};
\draw [fill=qqqqff] (4.735482645648355,-1.) circle (2.0pt);
\draw[color=qqqqff] (4.158074595203722,-0.6660098629189906) node {(5,6,10)};
\draw[color=ffqqqq] (-9.723956700951248,-3.1689371330180154) node {$P_{8}$};
\draw [fill=qqqqff] (-8.197747451591562,-3.) circle (2.0pt);
\draw[color=qqqqff] (-8.207931869893638,-3.460251708397712) node {(2,8,11)};
\end{scriptsize}
\end{tikzpicture}

\caption{A movable projective incidence theorem illustrating the geometric non-realizability of the arrangement $\mathcal{C}_{1}$.}
\label{fig:3}
\end{center}
\end{figure}
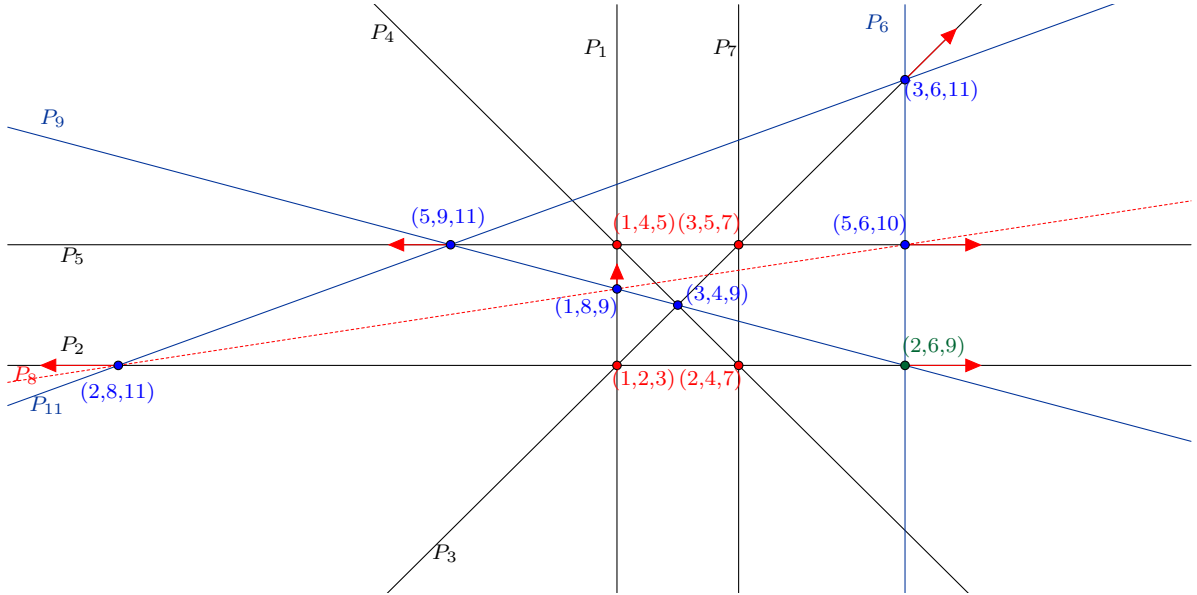

Let us now explain shortly why certain pseudoline arrangements are not realizable with straight lines. Recall that a \textbf{projective base} in the real projective plane is a set of four points in general position, i.e., no three points are collinear. We start with $\mathcal{C}_{1}$ and we use its corresponding tripels as points. When we pick a projective base of the points $(1,2,3), (1,4,5), (3,5,7)$, and  $(2,4,7)$, we have determined six lines $P_{1}, P_{2}, P_{3}, P_{4}, P_{5}$, and $P_{7}$ because two indices occur in former triples.
Now we can construct, as intersections,
the points $(1,6,7)$ and $(3,4,9)$. 
A movable point $(2,6,9)$ on line $P_{2}$ 
determines the points $(1,8,9), (5,9,11), (5,6,10)$,
and $(3,6,11)$. 
This defines in turn point $(2,8,11)$ and line $P_{8}$. 
However,
line $P_{8}$ passes through point $(5,6,10)$ because of a movable (!) projective incidence theorem, which is an obstruction. It is a {\it movable projective incidence theorem} because the incidence of line $P_{8}$ and point $(5,6,10)$
remains when you move point $(2,6,9)$ on line $P_{2}$. You can check this with any dynamic geometric software or you
provide an algebraic proof as follows. 
In our illustration in Figure \ref{fig:3}, we have chosen the projective base (four red points) as vertices of a square 
with edges of length 2 parallel to the coordinate axes and with midpoint as the origin. Thus the coordinates of all points are very easy and showing the final collinearity becomes an easy excercise:

$(1,2,3) = (-1,-1), \quad
(1,4,5) = (-1,1), \quad
(3,5,7) = (1,1), \quad
(2,4,7) = (1,-1), \quad
(3,4,9) = (0,0),$

$
(2,6,9) = (v,-1), v \neq 1, v \neq -1,v \neq 0, \quad (5,9,11)=(-v,1), \quad (1,8,9) = (-1, 1/v), \quad (5,6,10)=(v,1),
$

$(3,6,11)=(v,v), \quad (2,8,11)=(-v\cdot(3+v)/(v-1),-1).
$

$P_{11}: y = (v-1)x/2v + (v+1)/2, \quad P_{8}: y = (v-1)x/v(1+v) + 2/(1+v) 
$

We have not drawn $P_{10}$ since it does not play any role for finding the obstruction.

The second example $\mathcal{C}_{2}$ can be realized with straight lines. 
Here a movable projective incidence theorem works in favor of a realization. Pick again a projective base of the points $(1,2,3), (1,4,5), (3,5,7)$,  and  $(2,4,7)$ in a dynamical drawing software like Cinderella. After picking a movable point $(2,6,9)$ on line $P_{2}$, we see a movable projective incidence theorem when we forget about the two points $(7,9,10)$ and $(4,6,12)$. A corresponding line arrangement can be seen in Figure \ref{ex1}.

Let us point out here that $\mathcal{C}_{6}$ corresponds to B\"or\"oczky's arrangement of $12$ lines and $\mathcal{C}_{7}$ corresponds to the arrangement depicted in Figure \ref{movable}. As we can see, the arrangement $\mathcal{C}_{7}$ has $3$ axes of symmetry.

\section{New non-isomorphic line arrangements}   
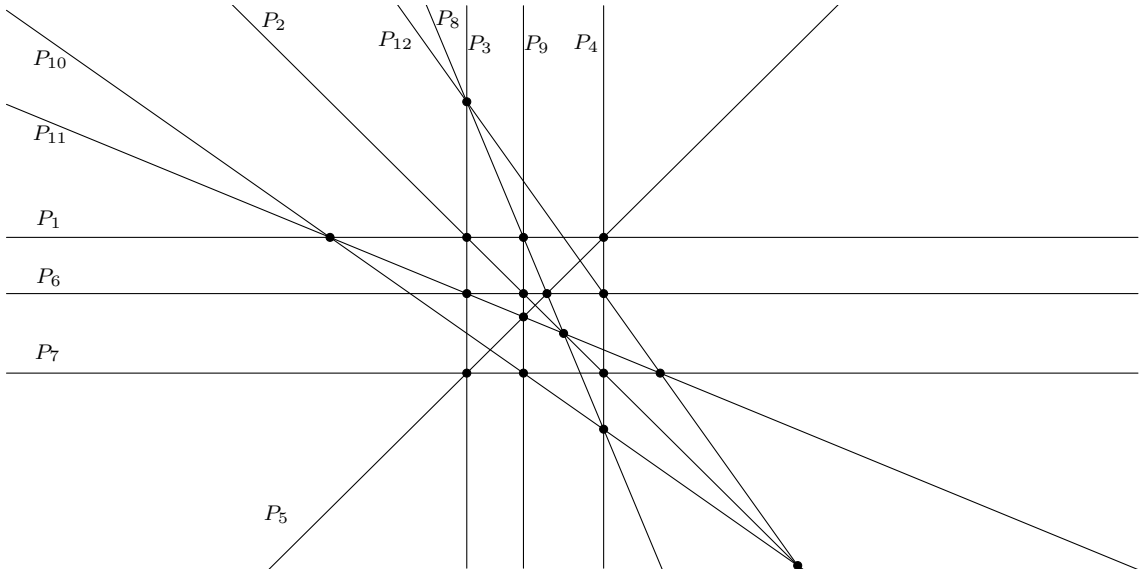
\begin{figure}[h]
\centering
\begin{tikzpicture}[line cap=round,line join=round,>=triangle 45,x=1.0cm,y=1.0cm,scale=0.6]
\clip(-11.08909649849381,-5.311877764308603) rectangle (13.712920118917216,7.1194508581101585);
\draw (-1.,-5.311877764308603) -- (-1.,7.1194508581101585);
\draw (2.,-5.311877764308603) -- (2.,7.1194508581101585);
\draw [domain=-11.08909649849381:13.712920118917216] plot(\x,{(-3.-0.*\x)/3.});
\draw [domain=-11.08909649849381:13.712920118917216] plot(\x,{(--6.-0.*\x)/3.});
\draw [domain=-11.08909649849381:13.712920118917216] plot(\x,{(-0.--3.*\x)/3.});
\draw [domain=-11.08909649849381:13.712920118917216] plot(\x,{(--3.-3.*\x)/3.});
\draw (0.24246561506204212,-5.311877764308603) -- (0.24246561506204212,7.1194508581101585);
\draw [domain=-11.08909649849381:13.712920118917216] plot(\x,{(--2.272603154813874-0.*\x)/3.});
\draw [domain=-11.08909649849381:13.712920118917216] plot(\x,{(--1.3313927293012882-1.242465615062042*\x)/0.5150687698759159});
\draw [domain=-11.08909649849381:13.712920118917216] plot(\x,{(--0.42614165563666995-0.5150687698759159*\x)/1.2424656150620421});
\draw [domain=-11.08909649849381:13.712920118917216] plot(\x,{(--10.75176620366536-4.239581524425743*\x)/3.});
\draw [domain=-11.08909649849381:13.712920118917216] plot(\x,{(-3.512184679239616-3.*\x)/4.239581524425742});
\begin{scriptsize}
\draw [fill=black] (-1.,2.) circle (2.5pt);
\draw [fill=black] (-1.,-1.) circle (2.5pt);
\draw[color=black] (-0.7119691301454923,6.270482074140097) node {$P_{3}$};
\draw [fill=black] (2.,2.) circle (2.5pt);
\draw [fill=black] (2.,-1.) circle (2.5pt);
\draw[color=black] (1.6226950257721715,6.2856422309967055) node {$P_{4}$};
\draw[color=black] (-10.18706716552562,-0.5515885113336136) node {$P_{7}$};
\draw[color=black] (-10.156746851812406,2.4046420757049938) node {$P_{1}$};
\draw[color=black] (-5.1690552459883055,-4.114225372636551) node {$P_{5}$};
\draw[color=black] (-5.229695873414738,6.770767250408169) node {$P_{2}$};
\draw [fill=black] (0.24246561506204212,0.7575343849379579) circle (2.5pt);
\draw[color=black] (0.5311637320963806,6.270482074140097) node {$P_{9}$};
\draw[color=black] (-10.156746851812406,1.1160287428932933) node {$P_{6}$};
\draw [fill=black] (0.757534384937958,0.757534384937958) circle (2.5pt);
\draw [fill=black] (0.24246561506204214,0.24246561506204214) circle (2.5pt);
\draw [fill=black] (0.24246561506204214,2.) circle (2.5pt);
\draw [fill=black] (-1.,0.757534384937958) circle (2.5pt);
\draw[color=black] (-1.40933634554947,6.8162477209779935) node {$P_{8}$};
\draw[color=black] (-10.12642653809919,4.284501525924417) node {$P_{11}$};
\draw [fill=black] (-3.9971159093637003,2.) circle (2.5pt);
\draw [fill=black] (-1.,4.997115909363701) circle (2.5pt);
\draw [fill=black] (0.24246561506204214,-1.) circle (2.5pt);
\draw [fill=black] (2.,0.757534384937958) circle (2.5pt);
\draw[color=black] (-2.561508266651694,6.361443015279747) node {$P_{12}$};
\draw[color=black] (-10.12642653809919,5.936958623294716) node {$P_{10}$};
\draw [fill=black] (2.,-2.239581524425742) circle (2.5pt);
\draw [fill=black] (3.2395815244257435,-1.) circle (2.5pt);
\draw [fill=black] (1.12225391796481,-0.12225391796481012) circle (2.5pt);
\draw [fill=black] (6.2535348026072715,-5.2535348026072715) circle (2.5pt);
\end{scriptsize}
\end{tikzpicture}
\caption{A geometric realization of the arrangement $\mathcal{C}_{2}$.}
\label{ex1}
\end{figure}
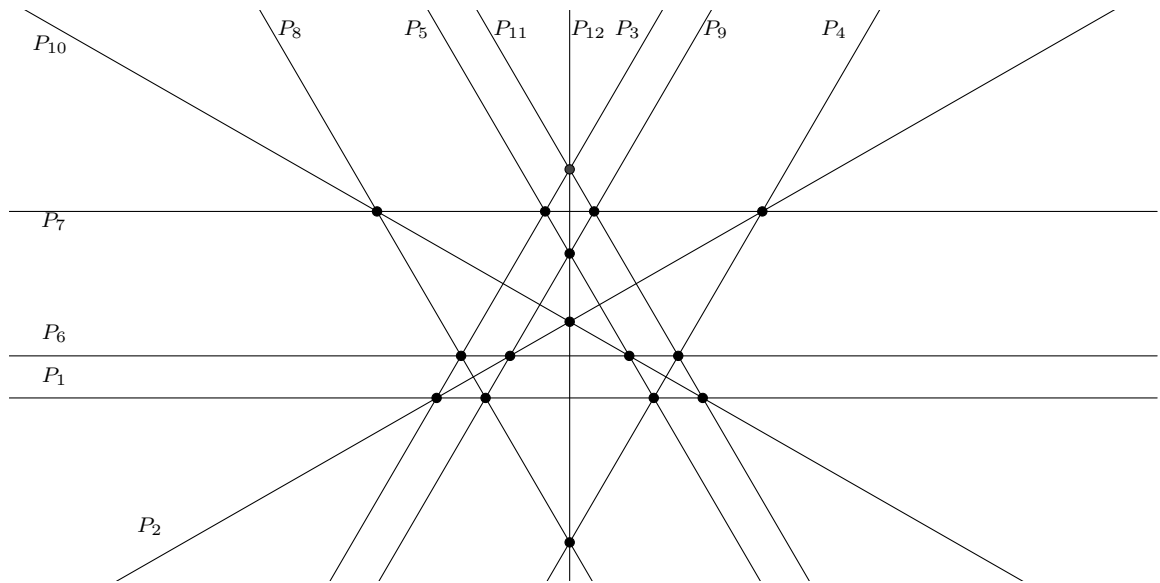
\begin{figure}[h]
\centering
\hskip 0.5 cm
\definecolor{uuuuuu}{rgb}{0.26666666666666666,0.26666666666666666,0.26666666666666666}
\begin{tikzpicture}[line cap=round,line join=round,>=triangle 45,x=1.0cm,y=1.0cm,scale=0.7]
\clip(-11.030042187083774,-4.463352965486819) rectangle (10.538621035893197,6.347346204953949);
\draw [domain=-11.030042187083774:10.538621035893197] plot(\x,{(--10.490381056766578--4.330127018922194*\x)/2.5});
\draw [domain=-11.030042187083774:10.538621035893197] plot(\x,{(--6.160254037844389-4.330127018922194*\x)/2.5});
\draw [domain=-11.030042187083774:10.538621035893197] plot(\x,{(--5.-0.*\x)/-5.});
\draw (-0.5,-4.463352965486819) -- (-0.5,6.347346204953949);
\draw [domain=-11.030042187083774:10.538621035893197] plot(\x,{(-3.169872981077803-2.5*\x)/-4.330127018922194});
\draw [domain=-11.030042187083774:10.538621035893197] plot(\x,{(--0.669872981077809-2.5*\x)/4.330127018922194});
\draw [domain=-11.030042187083774:10.538621035893197] plot(\x,{(-6.505652582187747-4.330127018922194*\x)/-2.5});
\draw [domain=-11.030042187083774:10.538621035893197] plot(\x,{(-2.175525563265561--4.330127018922194*\x)/-2.5});
\draw [domain=-11.030042187083774:10.538621035893197] plot(\x,{(--2.331118627451304-0.*\x)/0.9202336229782615});
\draw [domain=-11.030042187083774:10.538621035893197] plot(\x,{(-0.45469932893417664-0.*\x)/2.239299131065215});
\draw [domain=-11.030042187083774:10.538621035893197] plot(\x,{(--5.854893586094159-3.533181324006426*\x)/-2.0398831885108697});
\draw [domain=-11.030042187083774:10.538621035893197] plot(\x,{(-9.388074910100586-3.533181324006429*\x)/2.0398831885108697});
\begin{scriptsize}
\draw [fill=black] (-3.,-1.) circle (2.5pt);
\draw [fill=black] (2.,-1.) circle (2.5pt);
\draw [fill=uuuuuu] (-0.5,3.3301270189221936) circle (2.5pt);
\draw[color=black] (0.591459421140065,5.978200379621825) node {$P_{3}$};
\draw[color=black] (-1.5970479719003856,5.978200379621825) node {$P_{11}$};
\draw[color=black] (-10.179688410872274,-0.5873217994995183) node {$P_{1}$};
\draw[color=black] (-0.14683222952418334,5.978200379621825) node {$P_{12}$};
\draw[color=black] (-8.386694402116241,-3.43501816634733) node {$P_{2}$};
\draw[color=black] (-10.258791087729156,5.701341010622732) node {$P_{10}$};
\draw [fill=black] (-0.5,1.7362356290906618) circle (2.5pt);
\draw[color=black] (2.2394318556584767,5.978200379621825) node {$P_{9}$};
\draw[color=black] (-3.37685820118027,5.991384159097972) node {$P_{5}$};
\draw [fill=black] (-0.9601168114891293,2.533181324006428) circle (2.5pt);
\draw [fill=black] (-0.03988318851086779,2.5331813240064274) circle (2.5pt);
\draw[color=black] (-10.192872190348421,2.3394772442051766) node {$P_{7}$};
\draw [fill=black] (-4.119649565532607,2.533181324006429) circle (2.5pt);
\draw [fill=black] (3.119649565532608,2.533181324006426) circle (2.5pt);
\draw [fill=black] (-0.5,0.44337567297406427) circle (2.5pt);
\draw [fill=black] (-2.0797663770217376,-1.) circle (2.5pt);
\draw [fill=black] (1.0797663770217387,-1.) circle (2.5pt);
\draw [fill=black] (-1.6196495655326066,-0.20305430508423408) circle (2.5pt);
\draw [fill=black] (0.6196495655326084,-0.20305430508423422) circle (2.5pt);
\draw[color=black] (-10.179688410872274,0.23007252802161265) node {$P_{6}$};
\draw [fill=black] (-2.539883188510869,-0.203054305084234) circle (2.5pt);
\draw [fill=black] (1.5398831885108695,-0.20305430508423425) circle (2.5pt);
\draw[color=black] (4.46749058712737,5.991384159097972) node {$P_{4}$};
\draw[color=black] (-5.7499385068867825,5.978200379621825) node {$P_{8}$};
\draw [fill=black] (-0.5,-3.7362356290906575) circle (2.5pt);
\end{scriptsize}
\end{tikzpicture}
\caption{A geometric realization of the arrangement $\mathcal{C}_{7}$.}
\label{movable}
\end{figure}

\newpage

\section*{Acknowledgements}
We thank an anonymous referee for many suggestions that allowed us to improve the readability of our article.
The project was partially conducted when the second author was a fellow of SFB $45$ \textit{Periods, moduli spaces and arithmetic of algebraic varieties} and he was partially supported by National Science Centre Poland Grant 2014/15/N/ST1/02102. 
We wish to thank Torsten Hoge for useful remarks.

All figures in the paper are made with \emph{GeoGebra}.

%

\bigskip
   Piotr Pokora,
   Instytut Matematyki,
   Pedagogical University of Cracow,
   Podchor\c a\.zych 2,
   PL-30-084 Krak\'ow, Poland.

Current Address:
    Institut f\"ur Algebraische Geometrie,
    Leibniz Universit\"at Hannover,
    Welfengarten 1,
    D-30167 Hannover, Germany. \\
\nopagebreak
   \textit{E-mail address:} \texttt{piotrpkr@gmail.com}
   
\bigskip
J\"urgen Bokowski,
Department of Mathematics, Technische Universit\"at Darmstadt,
Schlossgartenstrasse~7, D-64289 Darmstadt, Germany.
\nopagebreak
   \textit{E-mail address:} \texttt{juergen.bokowski@gmail.com}


\end{document}